\documentclass[english,american]{amsart}
\usepackage[T1]{fontenc}
\usepackage[latin9]{inputenc}
\usepackage[letterpaper]{geometry}
\usepackage{float}
\usepackage{amsthm}
\usepackage{amstext}
\usepackage{amssymb}

\makeatletter

\providecommand{\tabularnewline}{\\}
\floatstyle{ruled}
\newfloat{algorithm}{tbp}{loa}
\floatname{algorithm}{Algorithm}

\numberwithin{equation}{section}
\numberwithin{figure}{section}
\theoremstyle{plain}
\newtheorem{thm}{Theorem}[section]
  \theoremstyle{definition}
  \newtheorem{problem}[thm]{Problem}
  \theoremstyle{remark}
  \newtheorem{rem}[thm]{Remark}
  \theoremstyle{plain}
  \newtheorem{conjecture}[thm]{Conjecture}

\usepackage{listings}

\makeatother

\usepackage{babel}

\makeatother

\usepackage{babel}

\begin{document}

\title[Critical branching random walk]{Critical branching random walk in an IID environment}

\author{János Engländer and Nándor Sieben}

\address{Department of Mathematics, University of Colorado, Boulder, CO-80309-0395,
and Department of Mathematics and Statistics, Northern Arizona University,
Flagstaff, AZ-86011.}

\email{Janos.Englander@Colorado.edu, nandor.sieben@nau.edu}

\urladdr{http://euclid.colorado.edu/~englandj/MyBoulderPage.html, http://jan.ucc.nau.edu/~ns46}



\thanks{The authors thank S. Kuznetsov and R. St. Laurent for helpful conversations.}

\date{\today}
\begin{abstract}
Using a high performance computer cluster, we run simulations regarding
an open problem about $d$-dimensional critical branching random walks
in a random IID environment The environment is given by the rule that
at every site independently, with probability $p>0$, there is a \emph{cookie},
completely suppressing the branching of any particle located there.

The simulations suggest \emph{self averaging}: the asymptotic survival
probability in $n$ steps is the same in the annealed and the quenched
case; it is $\frac{2}{qn}$, where $q:=1-p$. This particular asymptotics
indicates a non-trivial phenomenon: the tail of the survival probability
(both in the annealed and the quenched case) is the same as in the
case of non-spatial unit time critical branching, where the branching
rule is modified: branching only takes place with probability $q$
for every particle at every iteration. 
\end{abstract}
\maketitle
\tableofcontents{}

\section{Introduction}

\label{S:intro}

In \cite{E} a spatial branching model has been studied, where the
underlying motion is a $d$-dimensional ($d\ge1$) Brownian motion,
the particles perform dyadic branching, and the branching rate is
affected by a random collection of reproduction suppressing sets dubbed
\textit{mild obstacles}. In fact the obstacle configuration was given
by the union of balls with fixed radius, where the centers of the
balls form a Poisson point process. The radius $r$ plays no role
in the results, but the Poisson intensity $\nu>0$ does. The main
result of \cite{E} is the quenched Law of Large Numbers for the population
for all $d\ge1$. The environment $\omega$ (with law $P^{\omega}$)
has the property that the branching rate is $\beta_{1}$ inside the
obstacles and $\beta_{2}$ outside of them, with $0<\beta_{1}<\beta_{2}.$
The branching process given the environment $\omega$ was denoted
by $Z_{\cdot}(\omega)$ and the total population size at $t\ge0$
was denoted by $|Z_{t}(\omega)|$. Define the average growth rate
by $r_{t}=r_{t}(\omega):=\frac{\log|Z_{t}(\omega)|}{t}$. In \cite{E}
it was shown that for almost every environment, \[
\lim_{t\to\infty}(\log t)^{2/d}(r_{t}-\beta_{2})=-c(d,\nu)\quad\text{\ensuremath{P^{\omega}-}probability},\]
where $c(d,\nu)$>0 is an explicit constant. (This is a kind of LLN,
because the expectation of the total population size obeys the same
asymptotics.)

It has also been shown that the branching Brownian motion with mild
obstacles spreads less quickly than ordinary branching Brownian motion,
and an upper estimate for its radial speed has been provided.

More general offspring distributions (beyond the dyadic one considered
in the main theorems) were also discussed in \cite{E}. In particular,
the following question was posed. Consider the model where the offspring
distribution is \textit{critical}. One can easily prove (see Theorem
\ref{ext} below) that, despite of the obstacles, the system still
dies out with probability one.
\begin{problem}
What is the rate of decay for the survival probability? Is it still
of order $C/n$ as in the obstacle-free (non-spatial) case? 
\end{problem}
In the present paper we are going to investigate this problem in a
discretized setting. More precisely, we consider a modified version
of the model, by replacing the Poisson point process with IID probabilities
on $\mathbb{Z}^{d}$, as follows. Given an environment, the initial
particle, located at the origin, first moves according to a nearest
neighbor simple random walk, and immediately afterwards, the following
happens to her:
\begin{enumerate}
\item If there is no cookie at the new location, the particle either vanishes
or splits into two offspring particles, with equal probabilities.
\item If there is a cookie at the new location, nothing happens to the particle.
\end{enumerate}
The new generation then follows the same rule in the next unit time
interval and produces the third generation, etc.

Let $p\in[0,1]$. In the sequel $\mathbb{P}_{p}$ will denote the
law of the cookies and $P^{\omega}$ will denote the law of the BRW
given the environment $\omega$. So, if ${\mathbf{P}}_{p}$ denotes
the `mixed' law in the environment with cookie probability $p$, we
have \[
\mathbf{P}_{p}(\cdot)=\mathbb{E}_{p}P^{\omega}(\cdot).\]
 Following standard terminology, $P^{\omega}$ and $\mathbf{P}_{p}$
will be called \textit{quenched} and \textit{annealed} probabilities,
respectively.

We close this section with a remark which is essentially taken from
\cite{E} with slight adaptations.
\begin{rem}
An alternative view on our setting is as follows. Let $K$ denote
the (random) set of those lattice points where there is a cookie.
Then, our model can be viewed as a \textit{catalytic }BRW as well
--- the catalytic set is then  $K^{c}$ (in the sense that branching
is `made possible' there). Catalytic spatial branching (mostly for
superprocesses though) has been the subject of vigorous research in
the last twenty years initiated by Dawson, Fleischmann and others
--- see the survey papers \cite{Klenkesurvey} and \cite{DF} and
references therein. In those models the individual branching rates
of particles moving in space depend on the amount of contact between
the particle (`reactant') and a certain random medium called the catalyst.
The random medium is usually assumed to be a `thin' random set (that
could even be just one point) or another superprocess. Sometimes `mutually'
or even `cyclically' catalytic branching is considered \cite{DF}.

Our model is simpler than most catalytic models as our catalytic/blocking
areas are fixed, whereas in several catalytic models they are moving.
On the other hand, while for catalytic settings studied so far results
were mostly only qualitative we are aiming to get quite sharp \textit{quantitative}
results. 

For the discrete setting there are much less results available. One
example%
\footnote{A further example of the discrete setting is \cite{AB}.%
} is \cite{KestenSidoravicius}, where the branching particle system
on $\mathbb{Z}^{d}$ is so that its branching is catalyzed by another
autonomous particle system on $\mathbb{Z}^{d}$. There are two types
of particles, the $A$-particles (`catalyst') and the $B$-particles
(`reactant'). They move, branch and interact in the following way.
Let $N_{A}(x,s)$ and $N_{B}(x,s)$ denote the number of $A$- (resp.
$B$-)particles at $x\in\mathbb{Z}^{d}$ and at time $s\in[0,\infty)$.
(All $N_{A}(x,0)$ and $N_{B}(x,0)$, $x\in\mathbb{Z}^{d}$ are independent
Poisson variables with mean $\mu_{A}$ ($\mu_{B}$).) Every$A$-particle
($B$-particle) performs independently a continuous-time random walk
with jump rate $D_{A}$ ($D_{B}$). In addition a $B$-particle dies
at rate $\delta$, and, when present at $x$ at time $s$, it splits
into two in the next $ds$ time with probability $\beta N_{A}(x,s)ds+o(ds)$.
Conditionally on the system of the $A$-particles, the jumps, deaths
and splitting of the $B$-particles are independent. For large $\beta$
the existence of a critical $\delta$ is shown separating local extinction
regime from local survival regime. 

Finally, note that while the continuous equivalent of an IID trap
configuration on the lattice is a Poisson trap configuration on $\mathbb{R}^{d}$,
there is an important difference between the two. The discrete setting
has the advantage that the difference between the sets $K$ and $K^{c}$
is no longer relevant. Indeed, in the discrete case the complement
is also IID with a different parameter (self-duality), whereas in
the continuous setting this nice duality is lost as the `Swiss cheese'
$K^{c}$ is not the same type of geometric object as $K$; the latter
is the case, for example, in \cite{E}.
\end{rem}

\section{Preliminaries}

In this section we present two simple statements which are relatively
easy to verify rigorously. Let $S_{n}$ denote the event of \textit{survival}
for $n\ge0$. That is, $S_{n}=\{Z_{n}\ge1\}$, where $Z_{n}$ is the
population size at time $n$. 
\begin{thm}[Monotonicity]
\label{monotonicity} Let $0\le p<\widehat{p}\le1$ and fix $n\ge0$.
Then \[
{\mathbf{P}}_{p}(S_{n})\le{\mathbf{P}}_{\widehat{p}}(S_{n}).\]
 \end{thm}
\begin{proof}
First notice that it suffices to prove the following statement: 
\begin{quotation}
\textit{Assume that we are given an environment with some `red' cookies
and some additional `blue' cookies. Then the probability of $S_{n}$
with the additional cookies is larger than or equal to the probability
without them.}
\end{quotation}
Indeed, one can argue by coupling as follows. Let $q:=1-p,\ \delta:=\widehat{p}-p$.
First let us consider the cookies that are coming with IID probabilities
and parameter $p$. These will be the `red' cookies. Now with probability
$\delta/q$ at each site independently, add a blue cookie. Then the
probability for any given site, that there is at least one cookie
there is $p+\delta/q-p\delta/q=p+\delta=\widehat{p}.$ Now delete
those blue cookies where there was a red cookie too. This way, the
red cookies plus the additional blue cookies together correspond to
parameter $\widehat{p.}$

We are thus going to prove the statement in italics now, using an
argument due to S. Kuznetsov. Consider the generating functions of
no branching and critical branching: $\varphi_{1}(z)=z$ and $\varphi_{2}(z)=\frac{1}{2}(1+z^{2})$,
and note that $\varphi_{1}\le\varphi_{2}$ on $\mathbb{R}$. Fix an
environment and define \[
u(n,x,N):=P_{n,x}(S_{N}^{c}),\]
 that is, the probability that the population emanating from a single
particle which is at time $n$ is located in $x$, becomes extinct
at time $N$. Then, if the particle moves to the random location $\xi_{n+1}$,
one has

\begin{align*}
u(n,x,N) & =E\sum_{i=0}^{2}p_{i}(\xi_{n+1})\left[P_{n+1,\xi_{n+1}}(S_{N}^{c})\right]^{i}\\
 & =E\sum_{i=0}^{2}p_{i}(\xi_{n+1})\left[u(n+1,\xi_{n+1},N)\right]^{i}=E\varphi^{\xi_{n+1}}[u(n+1,\xi_{n+1},N)],\end{align*}
where $p_{i}(\xi_{n+1})$ is the probability%
\footnote{So either $p_{1}=1$ or $p_{0}=p_{2}=1/2$, according to whether there
is no cookie there or there is one.%
} of producing $i$ offspring at the location $\xi_{n+1}$, and $\varphi^{\xi_{n+1}}$
is either $\varphi_{1}$ or $\varphi_{2}$.

Now consider two environments: one with only red cookies and another
one where there are also some additional blue cookies, and let us
denote the corresponding functions by $u_{1}$ and $u_{2}$. We have
\[
u_{1}(n,x,N)=E\varphi_{1}^{\xi_{n+1}}[u_{1}(n+1,\xi_{n+1},N)]\]
 and \[
u_{2}(n,x,N)=E\varphi_{2}^{\xi_{n+1}}[u_{2}(n+1,\xi_{n+1},N)].\]
 Clearly $u_{1}(N,x,N)=u_{2}(N,x,N)=0$. Therefore, using that $\varphi_{1}\le\varphi_{2}$
and by (backward) induction, $u_{2}\ge u_{1}$ for all $n=0,1,...,N-1$.
In particular, \[
u_{1}(0,x,N)\le u_{2}(0,x,N),\]
 finishing the proof. 
\end{proof}
We now give a precise statement and a rigorous proof concerning the
result about eventual extinction mentioned briefly above. 
\begin{thm}[\emph{Extinction}]
\emph{}\label{ext} Let $0\le p<1$ and let $A$ denote the event
that the population survives forever. Then, for ${\mathbb{P}}_{p}$-almost
every environment, $P^{\omega}(A)=0.$ \end{thm}
\begin{proof}
Let again $Z_{n}$ denote the total population size at time $n$ for
$n\ge1$. Then $Z$ is a martingale with respect to the canonical
filtration $\{\mathcal{F}_{n};n\ge1\}$. To see this, note that just
like in the $p=0$ case, $E(Z_{n+1}-Z_{n}\mid\mathcal{F}_{n})=0$,
as the particles that do not branch (due to the presence of cookies)
do not contribute to the increment. Being a nonnegative martingale,
$Z$ converges a.s. to a limit $Z_{\infty}$, and of course $Z_{\infty}$
is nonnegative integer valued. We now show that for ${\mathbb{P}}_{p}$-almost
every environment, $P^{\omega}(Z_{\infty}=0)=1$. Introduce the events
\begin{itemize}
\item $A_{k}:=\{Z_{\infty}=k\}$ for $k\ge1$,
\item $B$: branching occurs at infinitely many times $0<\sigma_{1}<\sigma_{2}<...$
\end{itemize}
Clearly, $A=\cup_{k\ge1}A_{k}=\{Z_{\infty}\ge1\}$. We first show
that \begin{equation}
\text{for}\ \mathbb{P}_{p}-\text{a.e. environment},\ P^{\omega}(B^{c}A)=0.\label{io}\end{equation}
 Clearly, it is enough to show that $\mathbf{P}_{p}(B^{c}A)=0.$

Now, $B^{c}A\subset C$, where $C$ denotes the event that there exists
a first time $N$ such that for $n\ge N$, there is no branching and
particles survive and stay in the region of cookies. On $C$, one
can pick randomly a particle starting at $N$, and follow her path;
this path visits infinitely many points $P^{\omega}$-a.s. whatever
$\omega$ is%
\footnote{Because for every $\omega$, it is true $P^{\omega}$-a.s., that every
particle that does not branch, has a path that visits infinitely many
points%
}. Since this path is independent of $\omega$ and $p<1$, the $\mathbb{P}_{p}$-probability
that it contains a cookie at each of its sites is zero. Hence $\mathbf{P}_{p}(C)=0,$
and (\ref{io}) follows.

On the other hand, for each $k\ge1$, there is a $p_{k}<1$, such
that the probability that the population size remains unchanged (it
remains $k$) at $\sigma_{m}$ is not more than $p_{k}$ for every
$m\ge1$, uniformly in $\omega$. Thus, \[
P^{\omega}(BA_{k})=P^{\omega}(B\cap\{Z_{\sigma_{m}}=k\ \text{for all large enough}\ m\})=0,\]
 whatever $\omega$ is. Using this along with (\ref{io}), we have
that for ${\mathbb{P}}_{p}$-almost every $\omega$, $P^{\omega}(A_{k})=P^{\omega}(B^{c}A_{k})+P^{\omega}(BA_{k})=0,\ k\ge1$,
and so $P^{\omega}(A)=0.$ 
\end{proof}
\begin{algorithm}[t]
\lstset{language=C++,basicstyle=\small,identifierstyle=\slshape}
\lstset{morekeywords={string}}
\lstset{emph={cookie,simulate}, emphstyle=\sffamily}
\begin{lstlisting}[frame=no,numbers=left,numberstyle=\tiny,lineskip=-2pt,aboveskip=0pt,belowskip=0pt]
typedef struct {
  Tcell cell;          // location
  int iter;            // number of iterations
} Tparticle;
typedef vector < Tparticle > Tparticles;    

void simulate (
  int dim,             // dimension  
  double p,            // cookie probability 
  int maxiter,         // maximum number of iteration steps 
  int &iter            // longest survival
) {
 board.clear();                                  // remove cookies
 iter = 0;                                       // iteration counter
 Tparticles particles;                           // all the particles
 particles.reserve (maxiter/10);                  // reserve room
 Tparticle particle;
 particle.cell=Tcell(dim, 0);                    // cell at origin
 particle.iter = 0;
 particles.push_back (particle);                 // particle at center
 while (!particles.empty () && iter < maxiter) { // alive particles
   int coord = mtr.randInt (dim - 1);            // random coordinate
   int dir = 2 * mtr.randInt (1) - 1;            // -1 or +1
   particles.back ().cell[coord] += dir;         // move
   if (cookie (p, particles.back ().cell)) {     // there is a cookie 
     particles.back ().iter++;                   // survived iteration 
     iter = max (iter, particles.back ().iter);  // longest survival
   }
   else                                          // no cookie
     if (mtr2.randInt (1)) {                     // probability 0.5
       particles.back ().iter++;                 // survived iteration
       iter = max (iter, particles.back ().iter);// longest survival
       particles.push_back (particles.back ());  // duplicate particle
     }
     else                                        // die
       particles.pop_back ();                    // remove the particle
 }
}

\end{lstlisting}

\caption{\label{alg:code}The annealed simulation function of the code.}

\end{algorithm}

\section{Implementation}

The code for our simulations are written in C++ using the MPIqueue
parallel library \cite{MPIqueue}. We ran the code on 96 cores using
a computing cluster containing Quad-Core AMD Opteron(tm) 2350 CPU's.
We used an implementation \cite{MTR} of the Mersenne Twister \cite{Mersenne}
to generate random numbers. The total running time for the simulations
was several months.

\subsection{Annealed simulation}

Algorithm 1 shows the C++ function that runs a single annealed simulation.
We are essentially implementing a `depth-first search.' Below we give
a detailed description of the code. 
\begin{itemize}
\item line 1: We define a data type to store particles. 
\item line 2: The location of the particle is stored in the \textsl{cell}
field, that is a vector with the appropriate dimensions. 
\item line 3: The \textsl{iter} field stores the number of iterations survived
by the particle. 
\item line 5: We define a data type to store all the living particles. 
\item line 7: The simulation function takes three input variables and one
output variable. 
\item line 8: The dimension of the space is the first input. 
\item line 9: The probability of a cookie at any given location is the second
input. 
\item line 10: The maximum number of allowed iterations is the third input. 
\item line 11: The output of the function is the maximum number of iterations
any particle survived. 
\item line 13: We erase all the cookies from the board. Every run of the
simulation uses a new cookie placement. 
\item line 14: The initial value of the output must be zero. 
\item line 15: We define a variable to store all our alive particles. 
\item line 16: We reserve some space to store the particles. Making the
reserved space too small results in unnecessary reallocation of the
variable which degrades performance. On the other hand, reserving
too much space can be a problem too since different CPU's compete
with each other for RAM. 
\item lines 18--19: The initial particle starts at the origin before the
iterations start. 
\item line 20: At the beginning we only have the initial particle. 
\item line 21: We run the simulation while we have alive particles and none
of them stayed alive for the maximum allowed number of iterations. 
\item lines 22--23: We generate a random direction. 
\item line 24: We move the last of our alive particles in this random direction. 
\item line 25: We call the \textsf{cookie} function to check if there is
a cookie at the new location of the particle. The cookie function
checks in the global variable \textsl{board} if any particle already
visited this location and as a result we know already whether there
is a cookie there. If no particle visited this location before, then
the function uses the cookie probability to decide whether to place
a cookie there or leave the location empty. This information is then
stored for future visitors. 
\item line 26: If there is a cookie at the new location, then the particle
survived one more iteration, so we increment the \textsl{iter} variable. 
\item line 27: It is possible that this is the longest surviving particle
so far, so we update the output variable. 
\item line 29: If there is no cookie at the new location, then the particle
splits or dies. 
\item line 30: We generate a random number to decide what happens. 
\item lines 31--32: If the particle splits, then it survives so we update
information about the number of iterations. 
\item line 33: The particle splits, so we place a copy of it into our collection
of particles as the last particle. 
\item lines 35--36: If the particle dies, then we remove it from our collection
of particles. 
\end{itemize}
The rest of our code takes care of the parallelization, data collection
and the calculation of survival probabilities. The program splits
the available nodes into a boss node and several worker nodes. The
boss assigns simulation jobs to the workers. The workers call the
simulation function several times. The boss node collects the results
of these jobs and calculates the survival probabilities using all
the available simulation runs. More precisely, $\mathbf{P_{p}}(S_{n})$
is estimated as the number of simulation runs with longest survival
value not smaller than $n$ divided by the total number of runs.

\subsection{Quenched simulation}

The code for quenched simulation is essentially the same with only
minor modifications. In this version, line 13 of he simulation function
is missing, since we do not want to replace the board at every simulation.

The other change in the simulation function is at line 25. In the
annealed case, every worker node has a local version of the board
and the \textsf{cookie} function can create the board on the fly.
In the quenched case, the worker nodes need to use the same board,
so the \textsf{cookie} function cannot generate the board locally.
The new version of the cookie function still stores information about
the already visited locations. On the other hand, if a location is
not visited yet, then the worker node asks the boss node whether this
new location has a cookie. The boss node first checks whether the
location was visited by any other particle at any other worker node.
If the location was visited, then the boss already has a record of
this location. Otherwise, the boss node uses the cookie probability
to decide whether the location should have a cookie. Essentially,
the boss node has the ultimate information about the board, but the
worker nodes keep partial versions of the board and only consult the
boss node when it is necessary.
\begin{rem}
In the quenched case, note that if $\widehat{\rho}_{n}$ denotes the
relative frequency of survivals (up to $n$) after $r$ runs for a
fixed environment $\omega,$ that is,\[
\widehat{\rho}_{n}=\widehat{\rho}_{n}^{\omega}:=\frac{|\text{{survivals}}|}{r},\]
then using our method of simulation, the random variables $\widehat{\rho}_{n}$
and $\widehat{\rho}_{m}$ are not independent for $n\neq m$, because
the data are coming from the same $r$ runs.
\end{rem}
Similarly, in the annealed case, for a fixed environment and a fixed
run, the random variables $\mathbf{\mathbf{1}}_{S_{n}}$and $\mathbf{1}_{S_{m}}$
are not independent for $n\neq m$.

\section{Simulation results}

\begin{figure}
\input{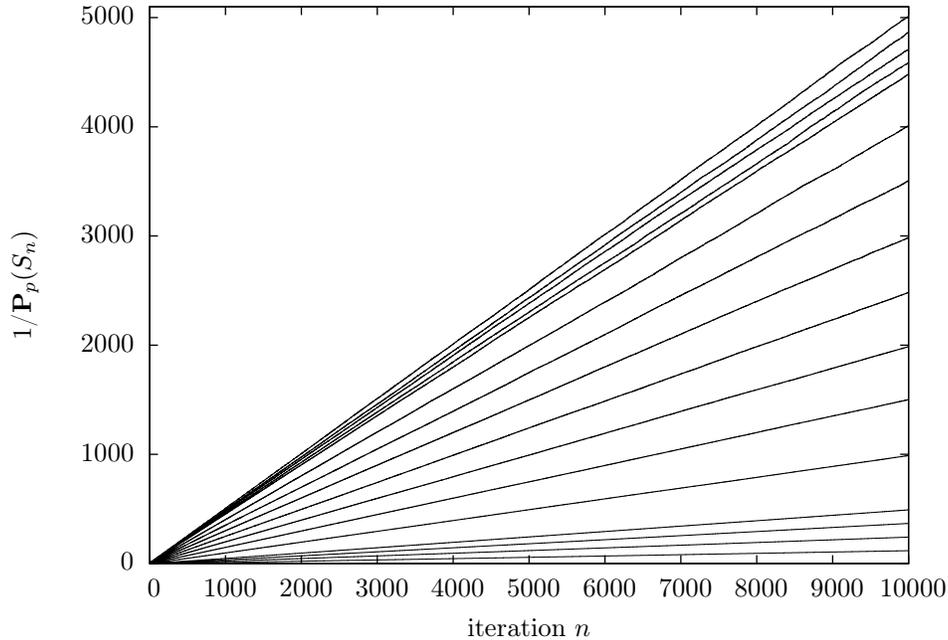}

\caption{\label{fig:annealed}Results for an annealed 2-dimensional simulation.
The graph shows the reciprocal of the survival probability as a function
of the number of iterations. Each line represents a different cookie
probability. One such line is the result of $10^{8}$ runs of the
simulation with a newly generated cookie landscape.}

\end{figure}

\begin{figure}
\input{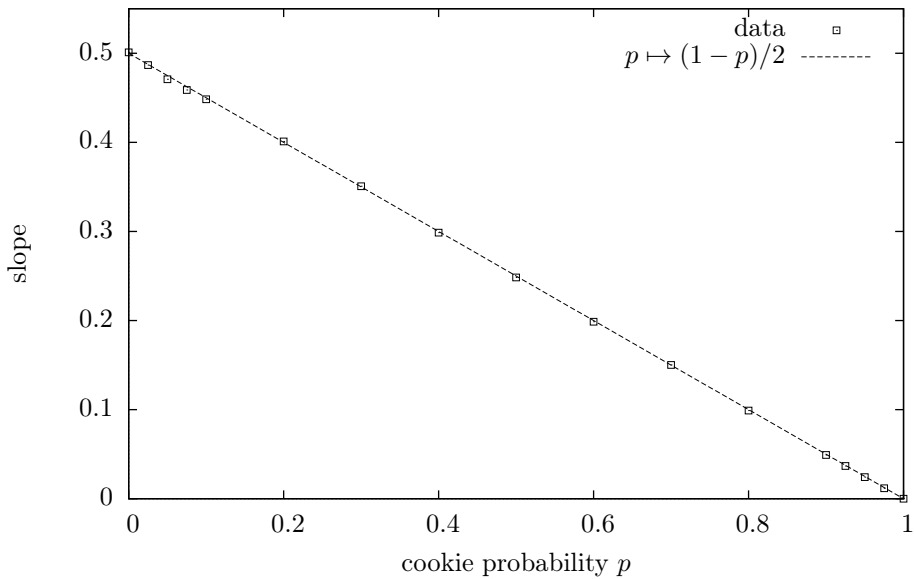}

\caption{\label{fig:annealed-1}Results for an annealed 2-dimensional simulation.
The graph shows the apparent slopes in Figure \ref{fig:annealed}
(i.e. the limits of the functions of Figure~\ref{fig:annealedPerT})
as a function of the cookie probability together with the graph of
$p\mapsto(1-p)/2$.}

\end{figure}

We ran our annealed simulation on $\mathbb{Z}^{d}$ with $d\in\{1,2,3\}$.
The 1-dimensional case turned out to be the most challenging. So we
start the description of our results with the 2-dimensional case.
The 3-dimensional case produced essentially the same output as the
2-dimensional case.

\subsection{Annealed simulation on $\mathbb{Z}^{2}$}

We executed $10^{8}$ runs allowing a maximum of $n_{\text{max}}=10000$
iterations with $p\in\{0,0.025,0.05,0.075,0.1,0.2,\ldots,0.9,0.925,0.95,0.975\}$.
For $p=1$ we used the known survival probability of $1$. Preliminary
results made it clear that the simulation is more sensitive for small
and large values of $p$. This is why we picked more of these values
instead of a uniformly placed set of values. The reciprocal of the
calculated survival probabilities are shown in Figure~\ref{fig:annealed}.
The figure suggests that $n\mapsto1/\mathbf{P}_{p}(S_{n})$ is asymptotically
linear. We calculated the slopes for these curves from the values
at $4n_{\text{max}}/5$ and $n_{\text{max}}$. These slope values
are shown in Figure~\ref{fig:annealed-1}.

To verify the correctness of our simulation we computed the exact
theoretical survival probabilities after the first two iterations.
It is easy to see that \foreignlanguage{english}{$\mathbf{P}_{p}(S_{1})=1/2+p/2$}
and\[
\mathbf{P}_{p}(S_{2})=3/8+11p/32+3p^{2}/16+3p^{3}/32.\]
Table~\ref{tab:Calculated-and-simulated-2} compares some of the
exact and simulated values.

\begin{table}
\begin{tabular}{|c|c|c|c|c|c|c|}
\hline 
$p$ & $0$ & $.5$ & $.975$ & $0$ & $.5$ & $.975$\tabularnewline
\hline 
$n$ & $1$ & $1$ & $1$ & $2$ & $2$ & $2$\tabularnewline
\hline 
exact  & $.5$ & $.75$ & $.9875$ & $.375$ & $.605469$ & $.975292$\tabularnewline
\hline 
simulated  & $0.50005$ & $.74998$ & $.98749$ & $.37501$ & $.605465$ & $.97528$\tabularnewline
\hline
\end{tabular}

\caption{\label{tab:Calculated-and-simulated-2}Exact and simulated survival
probability values $\mathbf{P}_{p}(S_{n})$ on $\mathbb{Z}^{2}$.
The simulated values are calculated from the data used in Figure~\ref{fig:annealed}. }

\end{table}

\subsection{Annealed simulation on $\mathbb{Z}^{1}$}

\begin{figure}
\input{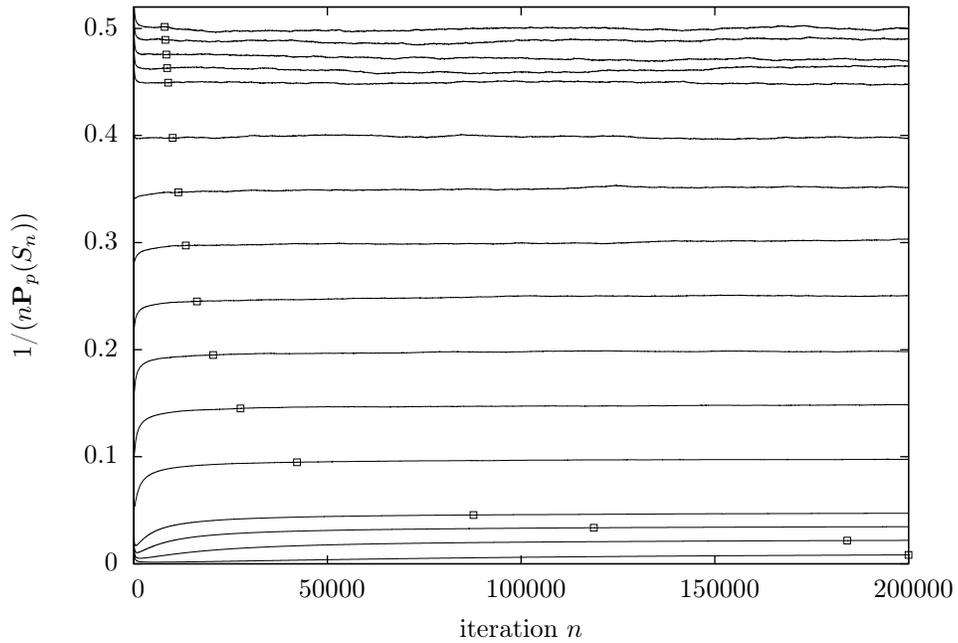}

\caption{\label{fig:avg1D-per-t-A}Annealed 1-dimensional simulation with $959965800$
runs. For small values of $p$ (lines at the top), a small iteration
number would actually give better results, because otherwise the survival
probability $\mathbf{P}_{p}(S_{n})$ becomes too small, even with
a huge number of runs. On the other hand, for large $p$ values (lines
at the bottom) one needs large iteration numbers because the convergence
is apparently slow. The squares represent the iteration thresholds
after which $\widehat{\rho}_{n}<\frac{1}{4000}$.}

\end{figure}

A 1-dimensional simulation with $10^{8}$ runs and $n_{\text{max}}=10000$
produces less satisfactory results as shown in Figure~\ref{fig:1Dannealed}.
The reasons behind this are explained in Subsection~\ref{sub:Interpretation.fluct}
below, where a discussion is given concerning the fluctuations of
the empirical curves in the figures. Essentially, in the annealed
case, small values of\foreignlanguage{english}{ $\mathbf{P}_{p}(S_{n})$}
result in large errors (see Subsection~\ref{sub:Interpretation.fluct}
for more explanation) and therefore we modified the original algorithm
by introducing a \textit{stopping rule: }when the estimated value
of \foreignlanguage{english}{$\mathbf{P}_{p}(S_{n})$} reaches a certain
small threshold value, we stop and do not simulate more iterations.
Fortunately, when larger threshold values needed, they are actually
large: we obtained slower convergence for large values of $p$, and,
clearly, for those values, the probability \foreignlanguage{english}{$\mathbf{P}_{p}(S_{n})$}
is large. The threshold value was set $1/4000$, based on trial and
error. This way, we stopped the iteration at $n_{\text{stop}}(p)$;
the slopes were then calculated from the values at $4n_{\text{stop}}(p)/5$
and $n_{\text{stop}}(p)$. See Figure \ref{fig:avg1D-per-t-A}.

After adjusting the algorithm by using the above stopping rule, the
curve indeed straightened out and the picture became very similar
to the 2-dimensional one in Figure~\ref{fig:annealed-1}.

\begin{figure}
\input{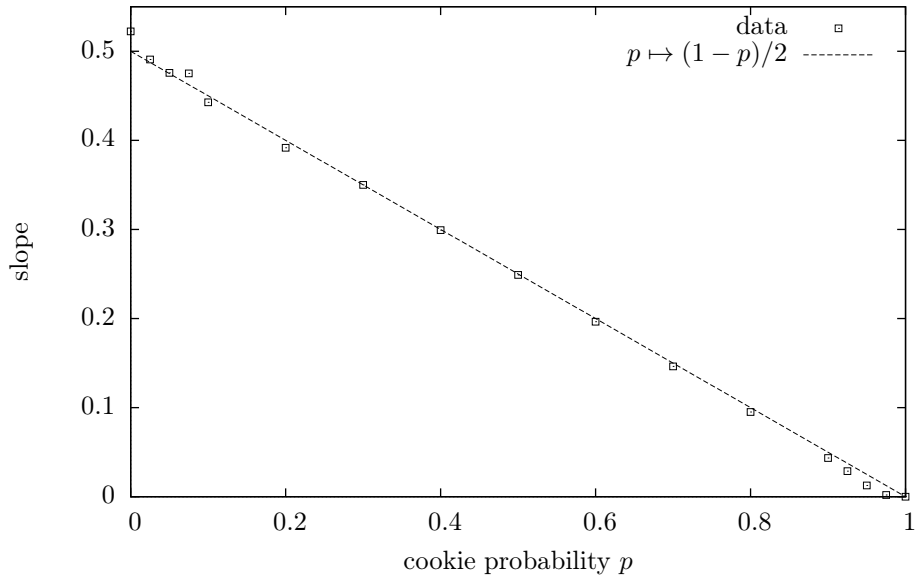}

\caption{\label{fig:1Dannealed}Results for an annealed 1-dimensional simulation.
Every parameter for this simulation is chosen to be the same as that
of Figure~\ref{fig:annealed-1} except the dimension.}

\end{figure}

To verify the correctness of our simulation we computed the exact
theoretical survival probabilities after the first two iterations.
It is easy to see that \foreignlanguage{english}{$\mathbf{P}_{p}(S_{1})=1/2+p/2$}
and\[
\mathbf{P}_{p}(S_{2})=3/8+5p/16+p^{2}/4+p^{3}/16.\]
Table~\ref{tab:Calculated-and-simulated} compares some of the exact
and simulated values.

\begin{table}
\begin{tabular}{|c|c|c|c|c|c|c|}
\hline 
$p$ & $0$ & $.5$ & $.975$ & $0$ & $.5$ & $.975$\tabularnewline
\hline 
$n$ & $1$ & $1$ & $1$ & $2$ & $2$ & $2$\tabularnewline
\hline 
exact  & $.5$ & $.75$ & $.9875$ & $.375$ & $.601563$ & $.975272$\tabularnewline
\hline 
simulated  & $0.50002$ & $.7499992$ & $.98749$ & $.37502$ & $.601569$ & $.975269$\tabularnewline
\hline
\end{tabular}

\caption{\label{tab:Calculated-and-simulated}Exact and simulated survival
probability values $\mathbf{P}_{p}(S_{n})$ on $\mathbb{Z}^{1}$.
The simulated values are calculated from the data used in Figure~\ref{fig:avg1D-per-t-A}. }

\end{table}

\subsection{Quenched simulation}

\begin{figure}
\input{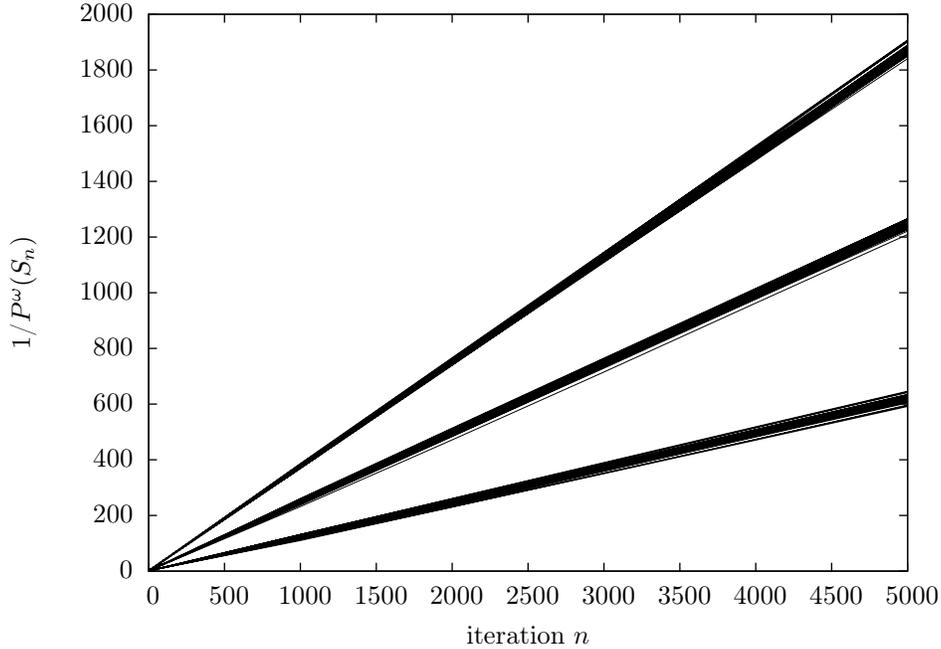}

\caption{\label{fig:quenched}Results for a quenched 2-dimensional simulation
with $10^{8}$ runs. Each line represents a different cookie landscape.
One such line is the result of $10^{8}$ runs of the simulation. The
lines are in three groups corresponding to three different cookie
probability. Each group has $50$ lines. The cookie probabilities
from top to bottom are $0.25$, $0.5$ and $0.75$. The total number
of simulations required for this graph is $3\cdot50\cdot10^{9}$.
The total running time was about $29$ hours.}

\end{figure}
 From the annealed simulation it has been clear that convergence is
much faster in two dimensions than in one dimension. Therefore, in
the quenched case we chose to present our results in the $d=2$ case.
In fact, we got qualitatively similar results for $d=1$.

In Figure~\ref{fig:quenched} we see three bundles corresponding
to three values of $p$. Those bundles are very thin, so essentially
the same thing happens for every realization, and the slopes of the
lines are roughly 3/8, 1/4 and 1/8 from top to bottom, corresponding
to $p=0.25$, $p=0.5$ and $p=0.75$, respectively. That is, for each
one of these values of $p$, the slope is the same as in the annealed
case.

Although Figure~\ref{fig:quenched} is about the $d=2$ case, we
have a similar simulation result for $d=1$; in fact \textit{we conjecture
that this qualitative behavior }(that is, the coincidence of the first
order asymptotics of the quenched and annealed survival probability)\textit{
will hold for all $d\ge1$. }

\section{Interpretation of the simulation results}

\subsection{Main finding}

Recall first the classic result due to Kolmogorov \cite[Formula 10.8]{H},
that for critical unit time branching with generating function $\varphi$,
as $n\to\infty$, \[
P(\text{survival up to}\ n)\sim\frac{2}{n\varphi''(1)}.\]
 As a particular case, let us consider now a non-spatial toy model
as follows. Suppose that branching occurs with probability $q\in(0,1)$,
and then it is critical binary, that is, consider the generating function
\[
\varphi(z)=(1-q)z+\frac{1}{2}q(1+z^{2}).\]
 It then follows that, as $n\to\infty$, \begin{equation}
P(\text{survival up to}\ n)\sim\frac{2}{qn}.\label{non-spatial}\end{equation}

Turning back to our spatial model, the simulations suggest (Figures~\ref{fig:annealed}
and \ref{fig:quenched}) the \textit{self averaging} property of the
model: as explained in the previous section, the asymptotics for the
annealed and the quenched case are the same. In fact, this asymptotics
is \textit{the same as the one in (\ref{non-spatial})}, where $p=1-q$
is the probability that a site has a cookie. In other words, despite
our model being spatial, in an asymptotic sense, the parameter $q$
simply plays the role of the branching probability of the above non-spatial
toy model. To put it yet another way, $q$ only introduces a `time-change.'

The intuitive picture behind this asymptotics is that \textit{there
is nothing either the environment or the BRW could do to increase
the chance of survival,} at least as far as the leading order term
is concerned (as opposed to well known models, for example when a
single Brownian motion is placed into random medium \cite{Sz}). Hence,
given any environment the particles move freely and experience branching
at $q$ proportion of the time elapsed (quenched case), and the asymptotics
agrees with the one obtained in the non-spatial setting as in (\ref{non-spatial}).
Furthermore, creating a `good environment' (annealed case) and staying
in the part of the lattice with cookies for very long would be `too
expensive.'

Note that whenever the total population size reduces to one, the probability
of that particle staying in the region of cookies is known to be much
less than $\mathcal{O}(1/n)$ (hard obstacle problem for random walk).
So the optimal strategy for this particle to survive is obviously
not to try to stay completely in that region and thus avoid branching.
Rather, survival will mostly be possible because of the potentially
large family tree stemming from that particle. In fact, the formula
$\mathbf{P}_{p}(S_{n})\sim\frac{2}{qn}$, together with the martingale
property of $|Z_{n}|$, implies that $\mathbf{E}_{p}(|Z_{n}|\mid S_{n})\sim\frac{q}{2}\cdot n.$ 

Notice that the straight lines on Figures \ref{fig:annealed-1} and
\ref{fig:1Dannealed} start at the value $1/2,$ that is, as $p\downarrow0$,
one gets the well known non-spatial asymptotics $2/n$, which is a
particular case of Kolmogorov's result above. We conclude that there
is apparently \textit{no discontinuity} at $p=0$ (no cookies) for
the quantity $\lim_{n\to\infty}nP(\text{survival up to}\ n)$.

\subsection{Interpretation of the fluctuations in the diagrams\label{sub:Interpretation.fluct}}

Since we estimated the reciprocal of the survival probabilities and
not the probabilities themselves both in the annealed and the quenched
case (Figures~\ref{fig:annealed} and \ref{fig:quenched}), one cannot
expect good approximation results when those probabilities are small.
Indeed, in the annealed case, if $\rho_{n}$:=$\mathbf{P}_{p}(S_{n})$
(with $p$ being fixed) and $\widehat{\rho}_{n}$ denotes the relative
frequency obtained from simulations, then the Law of Large Numbers
(LLN) only says, that if the number of runs is large, then the difference
$|\rho_{n}-\widehat{\rho}_{n}|$ is small. However, looking at the
difference of the reciprocals \[
\left|\frac{1}{\rho_{n}}-\frac{1}{\widehat{\rho}_{n}}\right|=\frac{|\rho_{n}-\widehat{\rho}_{n}|}{\rho_{n}\widehat{\rho}_{n}},\]
it is clear that a small $\rho_{n}$ value magnifies the error; in
fact the effect is squared as $\widehat{\rho}_{n}$ is close to $\rho_{n}$,
exactly because of the LLN. This effect is the reason of the `zigzagging'
of the line on Figure~\ref{fig:1Dannealed} for small $p$ values.
In fact, small $p$ values result in small $\rho_{n}$ values in light
of Theorem~\ref{monotonicity} and the continuity property mentioned
at the end of the previous subsection. Clearly, there is a competition
between $\rho_{n}$ being small (as a result of $p$ being small and
$n$ being large) on the one hand, and the large number of runs on
the other. The first makes the denominator small in $\frac{|\rho_{n}-\widehat{\rho}_{n}|}{\rho_{n}\widehat{\rho}_{n}},$
while the second makes the numerator small according to LLN.

Looking at Figure~\ref{fig:1Dannealed}, one notices another peculiarity
in the 1-dimensional setting. For large values of $p$, the empirical
curve is slightly under the straight line. The explanation for the
relatively poor fit is simply that the iteration number is not large
enough for the asymptotics to `kick in.' 

These arguments are bolstered by the experimental findings that increasing
the number of runs helps for small $p$ values, whereas increasing
the number of iterations helps for large ones. For example, in Figure~\ref{fig:avg1D-per-t-A}
we increased the maximal iteration number $n_{\text{max}}$ to $200000$
and plotted $n\mapsto(n\mathbf{P}_{p}(S_{n}))^{-1}$. One can see
that for small $p$ values it is beneficial to stop the iterations
earlier, but for large values large iteration numbers give better
results.

We do not have an explanation, however, for the deviation \emph{downward}
from the straight line (for large $p$ values) in Figure~\ref{fig:1Dannealed}.
We suggest, as an open problem, to find at least a heuristic explanation
for this phenomenon.

\begin{figure}
\input{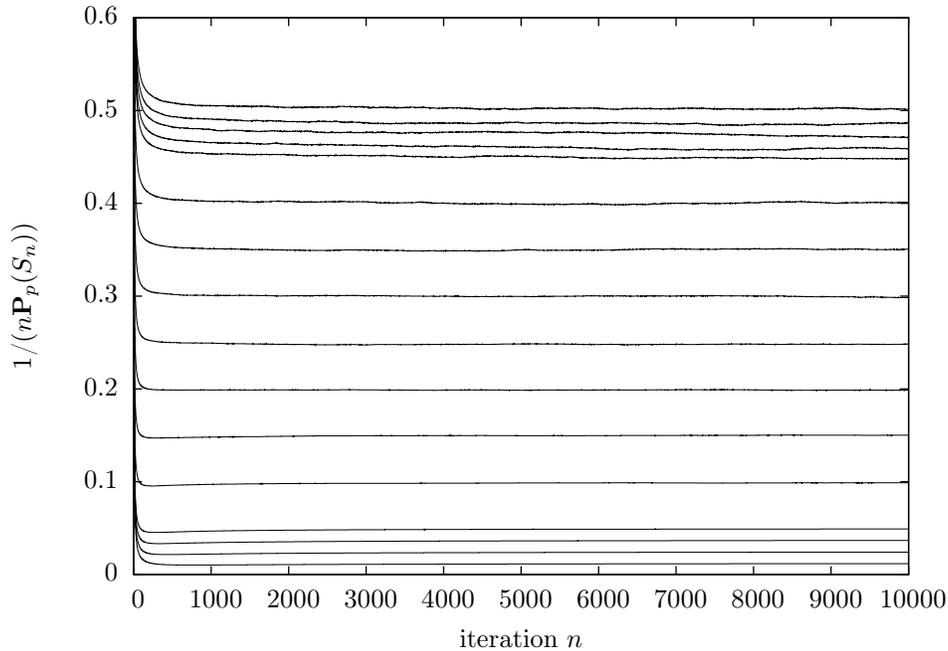}

\caption{\label{fig:annealedPerT}Results for an annealed 2-dimensional simulation
with $10^{8}$ runs. The graph shows the reciprocal of the survival
probability divided by the number of iterations as a function of the
number of iterations. The data used to create the graph is the same
as that of Figure~\ref{fig:annealed}. }

\end{figure}

Interestingly, for higher dimensions there is apparently a perfect
fit for large values of $p$, indicating that for higher dimensions
the asymptotics is much quicker then for $d=1$. Figure~\ref{fig:annealedPerT}
checks the assumption (for $d=2,$ annealed) that the reciprocal of
the survival probability is $\frac{qn}{2}+o(n)$ as $n\to\infty$,
by dividing the reciprocal of the survival probability by $n$. The
graphs convincingly show the existence of a limit, depending on $p$. 

\begin{figure}
\input{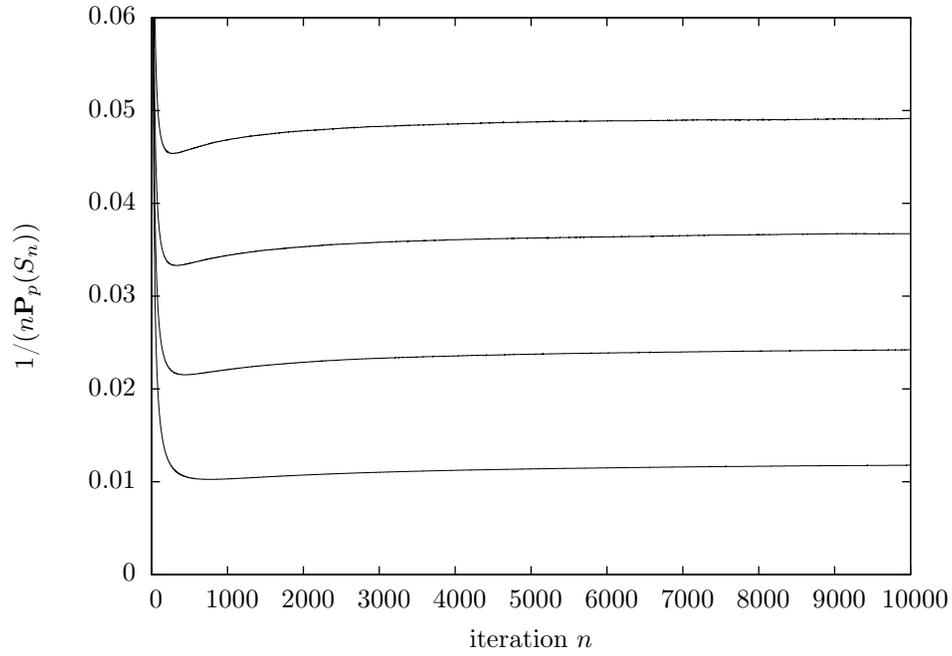}

\caption{\label{fig:ZoomannealedPerT}Zooming in at the bottom part of Figure
\ref{fig:annealedPerT} (i.e. large $p$ values): the convergence
is apparently from below.}

\end{figure}

\section{Beyond the first order asymptotics}

\begin{figure}
\input{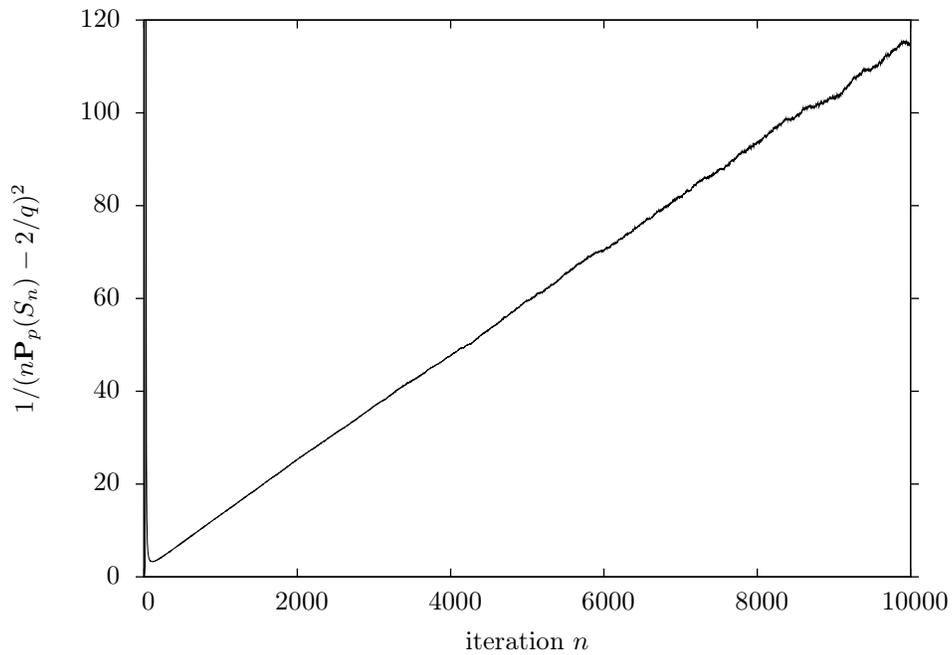}

\caption{\label{fig:1Dsquare}Annealed 1-dimensional simulation with 7,259,965,800
runs and $p=0.5$. }

\end{figure}

\begin{figure}
\input{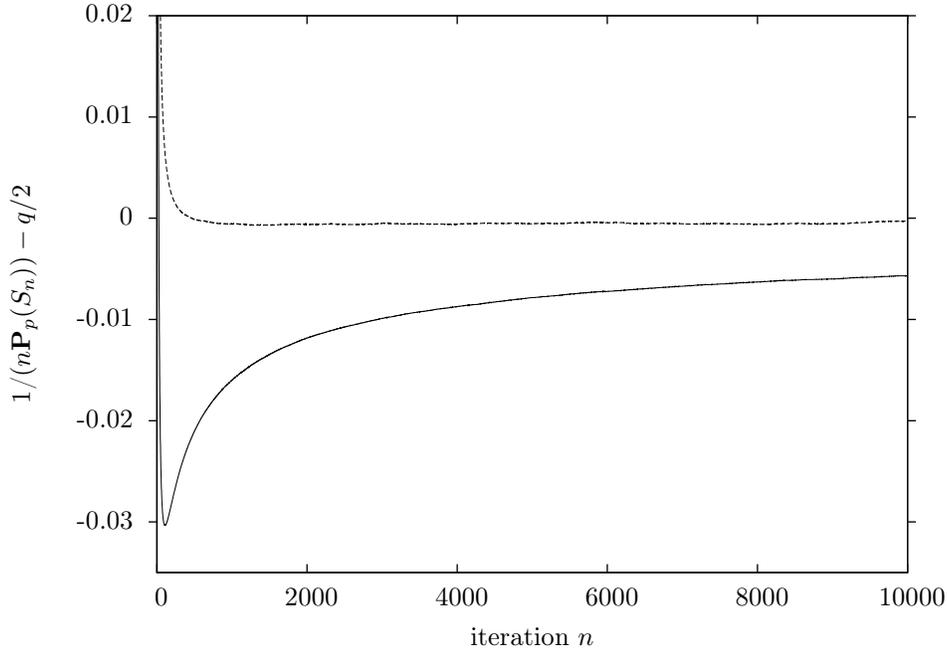}

\caption{\label{fig:Annealed 1D and D2 are different}Annealed simulation with
$p=0.5$. The solid curve shows the 1-dimensional result, the dashed
curve shows the 2-dimensional result.}

\end{figure}

In this section we will attempt to draw conclusions about more delicate
phenomena beyond the first order asymptotics, and those conclusions
will necessarily be less reliable than the ones in the previous sections.

\subsection{Two dimensions}

Consider again Figure~\ref{fig:annealedPerT}. Zooming in gives Figure~\ref{fig:ZoomannealedPerT}.
Looking at Figure \ref{fig:avg1D-per-t-A}, Figure~\ref{fig:annealedPerT}
and Figure~\ref{fig:ZoomannealedPerT}, for small $p$ values (top
lines) the convergence seems to be from above, and for large $p$
values it seems to be from below.

\subsection{One dimension}

For $d=1$ we obtained figures somewhat similar to the 2-dimensional
ones, which we summarize below without actually providing them. Simulation
seems to suggest that for not too small values of $p$, the convergence
is also from below; this is in line with the fact that, as we have
already discussed, on Figure~\ref{fig:1Dannealed} the 1-dimensional
empirical curve is \emph{below} the straight line for large $p$ values.
For very small $p$'s, the direction of the convergence is not clear
from the pictures. Although the convergence is apparently quicker,
the effects are `blurred' due to the magnification of error explained
earlier. 

The following conjecture is based on Figure~\ref{fig:1Dsquare}.
\begin{conjecture}[Second order asymptotics]
\label{con:[Second-order-asymptotics} The 1-dimensional annealed
survival probability obeys the following second order asymptotics:

\[
\mathbf{P}_{p}(S_{n})=\frac{2}{nq}+f(n),\]
where $\lim_{n\to\infty}f(n)\cdot n^{3/2}=C>0,$ and $C$ may depend
on $p$.
\end{conjecture}

\subsection{Comparison between one and two dimensions}

The annealed convergence to the limit $2/q$ seems to be quite different
for $d=1$ and $d=2$. Figure \ref{fig:Annealed 1D and D2 are different}
shows this difference, and in particular, it illustrates that in 1-dimension,
the convergence is slower and it is apparently from below for $p=0.5.$

\bibliographystyle{amsplain}
\bibliography{BRWsim}

\end{document}